\documentclass[centertags,12pt]{article}

\usepackage{amsmath}
\usepackage{amsthm}
\usepackage{amscd}
\usepackage{amssymb}
\usepackage{verbatim}

\newtheorem{theorem}{Theorem}
\newtheorem{proposition}[theorem]{Proposition}

\newtheorem{exmple}[theorem]{Example}

\newcommand{\Z}{{\bf Z}}
\renewcommand{\to}{\longrightarrow}
\newcommand{\Q}{{\bf Q}}

\DeclareMathOperator{\Aut}{Aut}
\DeclareMathOperator{\Comm}{Comm}
\DeclareMathOperator{\Out}{Out}
\DeclareMathOperator{\Inn}{Inn}

\DeclareMathOperator{\Mod}{Mod}

\DeclareMathOperator{\SL}{SL}
\DeclareMathOperator{\GL}{GL}
\DeclareMathOperator{\Sp}{Sp}
\DeclareMathOperator{\SP}{Sp}

\DeclareMathOperator{\T}{{\mathcal I}}

\DeclareMathOperator{\genus}{genus}

\DeclareMathOperator{\TG}{{\mathcal{TG}}}
\DeclareMathOperator{\TGS}{{\mathcal C}_{sep}}
\DeclareMathOperator{\CC}{{\mathcal{CC}}}
\newcommand{\mme}{{\rm Mod}^{\pm}}
\newcommand{\spe}{{\rm Sp}^{\pm}}

\title{The Torelli geometry and its applications
\\{\small RESEARCH ANNOUNCEMENT}}
\author{Benson Farb \thanks{The first author was 
supported in part by NSF grants DMS-9704640 and DMS-0244542. He would
also like to thank the Institut de Math\'{e}matiques de Bourgogne for
its hospitality and support during the writing of this paper.} \ and
Nikolai V. Ivanov \thanks{This work was done during the visit of the
second author to the University of Chicago. He would like to thank the
Department of Mathematics of the University of Chicago and the first
author for the hospitality.}}
\date{}

\begin{document}
\maketitle

Let $S$ be a closed orientable surface of genus $g$.
The {\em mapping class group\/} $\Mod (S)$ of $S$ is defined
as the group of isotopy classes of orientation-preserving diffeomorphisms $S\rightarrow S$.
We will need also the {\em extended mapping class group\/} $\mme (S)$
of $S$ which is defined as the group of isotopy classes of {\em all\/} diffeomorphisms
$S\rightarrow S$. Let us fix an orientation of $S$. Then
the algebraic intersection number provides a nondegenerate,
skew-symmetric, bilinear form on $H=H_1(S,\Z)$, called the {\em
intersection form.} The natural action of $\Mod (S)$ on $H$ preserves 
the intersection form. If we fix a symplectic basis in $H$, then
we can identify the group of symplectic automorphisms of $H$
with the integral symplectic group $\Sp(2g,\Z)$ and the
action of $\Mod (S)$ on $H$ leads to a natural surjective homomorphism
$$\Mod (S)\longrightarrow \Sp(2g,\Z).$$ 
The {\em Torelli group} of $S$, denoted by $\T_S$, is defined
as the kernel of this homomorphism; that is, $\T_S$ is the subgroup
of $\Mod (S)$ consisting of elements which act trivially on
$H_1(S,\Z)$. In particular, we have the following well-known exact sequence
\begin{eqnarray}
\label{equation:basic}
1\longrightarrow \T_S\longrightarrow 
\Mod (S) \longrightarrow \Sp(2g,\Z)\longrightarrow 1.
\end{eqnarray}

The group $\T_S$ plays an important role both in algebraic geometry and
in low-dimensional topology. At the same time most of the basic questions about $\T_S$
remain open. See, e.g., \cite{Jo1} and \cite{Ha1} for surveys. For
example, while $\T_S$ is finitely generated if the genus $g\geq 3$ by a famous
theorem of D. Johnson \cite{Jo2}, it is not known whether or not $\T_S$ 
has a finite presentation if $g\geq 3$. 

In this note we introduce a new geometric object related to the Torelli
group $\T_S$, which we call the {\em Torelli geometry} $\TG(S)$ of $S$
and announce several results about $\TG (S)$ and $\T_S$. Our main result
about $\TG (S)$ gives a complete description of its automorphisms (they
are all induced by diffeomorphisms of $S$). We also give a purely
algebraic characterization of some geometrically defined elements of
$\T_S$, namely, of the (powers of the) so-called {\em bounding twists}
and {\em bounding pairs}, as also of some geometrically defined
collections of such elements. When combined with the description of the
automorphisms of $\TG (S)$, these characterizations lead to several
results about the Torelli group $\T_S$, among which are the following.

\begin{itemize}
\item The automorphism group of $\T(S)$ is isomorphic to the extended mapping
class group $\mme (S)$; the outer automorphism group $\Out(\T(S))$ contains
the integral symplectic group $\SP(2g,\Z)$ as a subgroup of index $2$.  

\item The abstract commensurator $\Comm(\T(S))$ is isomorphic to the extended mapping class group 
$\mme (S)$.

\item The Torelli group is not arithmetic.
\end{itemize}

The last result is known. Actually, any normal subgroup of $\Mod (S)$
was proved to be non-arithmetic in \cite{i-ar}; for more general
results, proving a conjecture of the second author, see
\cite{farb-m}. However, we feel that our proof sheds a new light on this
result.

The details of the proofs will appear in \cite{FI}.  

\bigskip
\noindent
{\bf Torelli Geometry.} A recurring theme in the theories of both finite
and infinite groups is the investigation of a group by using its action
on an appropriate {\em geometry}. We will use the term {\em geometry} in
a narrow sense of a simplicial complex with some additional structure,
namely, with some vertices and simplices marked by some colors. For
example, the main construction of Bruhat-Tits theory applied to
$\SL_n(\Q_p), n>2$ gives an $(n-1)$-dimensional simplicial complex
$X_{n,p}$ (with marked vertices) whose group of simplicial automorphisms
preserving the marking is isomorphic to $\SL_n(\Q_p)$. Similar
constructions occur in the theory of finite simple groups.

For the mapping class group $\Mod(S)$ of a closed surface $S$, a useful
geometry is the {\em complex of curves}, denoted $\CC(S)$. This
simplicial complex, introduced by W. Harvey \cite{ha2}, is defined in
the following way: $\CC(S)$ has one vertex for each isotopy class of
simple (i.e. embedded) closed curve on $S$; further, $k+1$ distinct
vertices of $\CC(S)$ form a $k$-simplex if the corresponding isotopy
classes can be represented by disjoint curves. It is easy to see that
$\CC(S)$ has dimension $3g-4$; in contrast with the classical situations
it is not locally finite.

The complex $\CC(S)$ is a fundamental tool in the low-dimensional
topology.  The following theorem, which strongly supports the claim of
$\CC (S)$ for being the right geometry for $\Mod (S)$, is of special
importance for us.

\begin{theorem}[Ivanov \cite{Iv2}]
\label{theorem:ivanov}
Let $S$ be closed and let $\genus(S) \geq 2$. 
Then the group of simplicial automorphisms of $\CC(S)$ is isomorphic to
the extended mapping class group $\mme (S)$.
\end{theorem}

Our definition of a geometry for the Torelli group is guided by the
definition of the complex of curves $\CC (S)$. In order to arrive at a
complex intimately related to the Torelli group, it is natural to use as
vertices the isotopy classes of {\em separating curves} and of {\em
bounding pairs}, since they correspond to the natural generators of the
Torelli group.  By a {\em separating curve} $\gamma\subset S$ we mean a
homologically trivial curve, or equivalently, a curve which bounds a
subsurface in $S$.  A {\em bounding pair} in $S$ is a pair of disjoint,
non-isotopic, nonseparating curves in $S$ whose union bounds a
subsurface in $S$. It turns out that it is useful also to introduce a
marking.

\bigskip
\noindent
{\bf Definition (Torelli geometry).}  The {\em Torelli geometry} of $S$,
denoted $\TG(S)$, is the following simplicial complex with additional
structure:

\bigskip
\noindent
{\bf Vertices of \boldmath$\TG(S)$: }The vertices of $\TG(S)$ consist of:
\begin{enumerate}
\item Isotopy classes of separating curves in $S$, and
\item Isotopy classes of bounding pairs in $S$. 
\end{enumerate}

\bigskip
\noindent
{\bf Simplices of \boldmath$\TG(S)$: }A collection of $k\geq 2$ 
vertices in $\TG(S)$ forms a $(k-1)$-simplex if these vertices have
representatives which are 
mutually non-isotopic and disjoint.

\bigskip

We also endow $\TG(S)$ with a {\em marking}, which consists of the
following two pieces of data:

\bigskip
\noindent
{\bf Marking for \boldmath$\TG(S)$: }
\begin{enumerate}

\item  Each vertex in $\TG(S)$ is marked by its type: 
separating curve or bounding
pair.

\item A $2$-simplex $\Delta$ in $\TG(S)$ is marked if there are $3$
non-isotopic, disjoint, nonseparating curves 
$\gamma_1,\gamma_2,\gamma_3$ so that the vertices of $\Delta$ are
bounding pairs $\gamma_i \cup \gamma_{i+1}$ for $i=1,2,3\ ({\rm mod}\,
3)$.  
\end{enumerate}
\bigskip


It is easy to see that the group $\Mod(S)$ acts on $\TG(S)$ by
automorphisms, by which we mean simplicial automorphisms 
preserving the marking. Our main result about $\TG(S)$ is the following.

\begin{theorem}[Automorphisms of \boldmath$\TG(S)$]
\label{automorphisms of geometry} 
Suppose $\genus(S)\geq 5$. Then every 
automorphism of the Torelli geometry $T(S)$
is induced by a diffeomorphism of $S$.  
Further, the natural map
$$\mme (S)\longrightarrow \Aut(\TG(S))$$
is an isomorphism. 
\end{theorem}

The fact that $\Aut(\TG(S))$ is the extended mapping class group 
and not the Torelli group itself reflects intrinsic extra symmetries of
the Torelli group; indeed, exactly by this reason Theorem
\ref{automorphisms of geometry} will allow us to
compute the automorphisms group of $\TG (S)$ and the abstract
commensurator of $\TG (S)$.

\bigskip
\noindent
{\bf Connectedness theorems. }A key ingredient in the proof of Theorem
\ref{automorphisms of geometry} is the following result, which we believe is of
independent interest.

\begin{theorem}[\boldmath$\TG(S)$ is connected]
\label{connect1}
For a closed orientable surface $S$ of $\genus(S)\geq 3$ the complex
$\TG(S)$ is connected.
\end{theorem}

There is a basic and useful subcomplex of the complex of curves, 
which is in fact a subcomplex of $\TG(S)$ for closed $S$.

\bigskip
\noindent
{\bf Definition (\boldmath$\TGS(R)$).} Let $R$ be any compact orientable surface,
perhaps with nonempty boundary. We define $\TGS(R)$ to be the
simplicial complex whose vertices consist of isotopy classes of curves
in $R$ which bound a subsurface of positive genus with one boundary
component, and with $k+1$ distinct vertices forming a $k$-simplex if
they can be represented by disjoint curves on $S$.
\bigskip

Note that, for a closed surface $S$, the complex $\TGS(S)$ is just the
subcomplex of $\CC(S)$ spanned by the (isotopy classes of) separating curves.

\begin{theorem}[\boldmath$\TGS(R)$ is connected]
\label{connect2}
Suppose that $R$ is a compact surface with $\genus(R)\geq 3$. Then 
$\TGS(R)$ is connected.
\end{theorem}

Note that if $\gamma\subset R$ is a separating curve, and $\genus(R)\geq
3$, then $R\setminus \gamma$ has one component of genus at least two; in
particular there is a separating curve $\beta \in R$ disjoint from
$\gamma$ and bounding a subsurface of genus one. We will call such
separating curves the {\em genus one separating curves}. Hence 
Theorem \ref{connect2} follows from the following more strong result
which we prove: 

\bigskip
\noindent
{\em For $\genus(R)\geq 3$,the full subcomplex of $\TGS(R)$
spanned by the (isotopy classes of) genus one separating curves is connected.}  
\bigskip

Notice that both connectedness theorems do not hold if we omit the
restriction on the genus; in fact, if the genus is equal to $2$, then 
both complexes have an infinite number of vertices and no edges.

\bigskip
\noindent
{\bf The idea of proof of Theorem \ref{automorphisms of geometry}. } The main
step in the proof of Theorem \ref{automorphisms of geometry} is to show
that every automorphism of the Torelli geometry $\TG(S)$ canonically
induces an automorphism of the complex of curves $\CC(S)$ when
$\genus(S)\geq 5$.  Once this is proved, we can quote Ivanov's
Theorem (Theorem \ref{theorem:ivanov} above) and complete the proof.

The problem is to encode a nonseparating curve in $\CC(S)$ in terms of 
homologically trivial curves. The main difficulty results from the ``infiniteness'' of
this problem; for example, any nonseparating curve can be completed to a bounding pair in infinitely many different ways, none of which is more natural than any other. The key idea
is to consider all possible completions at once. 

Consider a vertex $\gamma$ of $\TG(S)$ represented by a bounding pair of
nonseparating curves $C_0,C_1$. We would like to single out one of
these curves using only the information available in $\TG(S)$. This can be done
by using a separating curve $D$ which has $i(D,C_0)\neq 0$ and
$i(D,C_1)=0$.  The isotopy class $\langle C_0\rangle$ obviously can be
recovered from the pair $(\gamma, \delta)$, where $\delta=\langle
D\rangle$. Therefore the isotopy class $\langle C_0\rangle$ can be
encoded by $(\gamma, \delta)$. For technical reasons, we use only
pairs which we will call {\em admissible} (see below).  Of course, 
we need to be able to tell when two such pairs encode the same isotopy class. 

\bigskip
\noindent
{\bf Definition (Admissible pair). }
Let $\gamma$ be a bounding pair  and let $\delta$ be a separating
curve. We call the pair $(\gamma, \delta)$ an 
{\em admissible pair} if both:
\begin{enumerate}
\item $\gamma$ is a vertex of some marked triangle
$\{ \gamma , \gamma' ,\beta\}$ of $T(S)$, and 
\item $\delta$ is connected by an edge of $T(S)$ with $\beta$ and is {\em not}
connected by an edge with either $\gamma$ or $\gamma'$. 
\end{enumerate}


We will say that the triangle $\{ \gamma , \gamma' ,\beta\}$ {\em
certifies the admissibility} of the pair $(\gamma, \delta)$. 

\bigskip
Let the marked triangle $\{ \gamma , \gamma' ,\beta\}$ be determined by
three curves $C_0$, $C_1$, $C_2$ (as in the definition of the marked
triangles), and suppose that $\beta$ is the isotopy class of the pair
$C_1$, $C_2$ (and, hence, $\gamma$, $\gamma'$ are the isotopy classes of
two pairs including $C_0$). Then $\delta$ can be represented by a curve
$D$ disjoint from both $C_1$ and $C_2$. Since $\delta$ is not connected
by an edge with $\gamma$, we have $i(D,C_0)\neq 0$. It follows, in
particular, that the isotopy class $\langle C_0\rangle$ of the
nonseparating curve $C_0$ is uniquely determined by the pair $(\gamma,
\delta)$.  We will say that the isotopy class $\langle C_0\rangle$ is
{\em encoded\/} by $(\gamma,
\delta)$. Obviously, such an encoding of $\langle C_0\rangle$ is far
from being unique. We will account for this non-uniqueness by using the
following {\em moves}.

\bigskip
\noindent
{\bf Two moves on admissible pairs: }
We consider the following two types of moves for admissible 
pairs $(\gamma,\delta)$.

\bigskip
\noindent
{\bf Type I move:} If a marked triangle 
$\{ \gamma , \gamma' ,\beta\}$ certifies the admissibility
of the pair $(\gamma, \delta)$, then replace $(\gamma, \delta)$
by $(\gamma', \delta)$.

\medskip
\noindent
{\bf Type II move:} If a marked triangle $\{ \gamma , \gamma' ,\beta\}$ 
certifies the admissibility of the pair $(\gamma, \delta)$ and simultaneously
certifies the admissibility of the pair $(\gamma, \delta')$,
then replace $(\gamma, \delta)$ by $(\gamma, \delta')$. 
\bigskip

The key result about the encodings of nonseparating curves by homologically
trivial ones is the following.

\begin{theorem}[Equivalence of encodings]
\label{encoding}
Two admissible pairs encode the same nonseparating curve if 
and only if they are connected by a sequence of moves of types I and II. 
\end{theorem}

This theorem is deduced from the above connectedness theorems, namely 
from Theorems
\ref{connect1}, \ref{connect2}.

Once this is proven, the proof that an automorphism of $\TG(S)$ induces
an automorphism of $\CC(S)$ follows easily, noting that the definition
of admissible pairs and of moves of types I and II 
are given entirely in terms of the Torelli geometry.

\bigskip
\noindent
{\bf Applications to the Torelli group. }In order to state the main results about the
Torelli groups, we need a counterpart of the extended mapping class group $\mme (S)$ 
for the symplectic groups. Namely, we need 
{\em the extended symplectic group} $\spe (2g, \Z)$, which is defined
as the group of automorphisms of $\Z^{2g}$ preserving the standard symplectic
form up to an overall sign. Clearly, $\spe (2g, \Z)$ contains 
$\Sp(2g,\Z)$ as a subgroup of index $2$.

Let $S$ be a closed surface of genus $g$. The conjugation action of
$\mme (S)$ on the normal subgroup $\T(S)$ induces homomorphisms $$\mme(S)\to
\Aut(\T(S)) \ \ \ \mbox{and}\ \ \
\spe (2g,\Z)\to \Out(\T(S))$$
where $\Aut(\T(S))$ and $\Out(\T(S))$ are the groups of automorphisms and
outer automorphisms, respectively, of $\T(S)$.

Our first application of
Theorem \ref{automorphisms of geometry} is the following.

\begin{theorem}[Automorphisms of \boldmath$\T(S)$]
\label{theorem:autos}
Let $S$ be a compact surface with $\genus(S)\geq 5$.  Then 
the natural homomorphisms 
$$\mme(S)\to \Aut(\T(S)) \ \ \ \mbox{and}\ \ \
\spe (2g,\Z)\to \Out(\T(S))$$
are isomorphisms. In fact there is an isomorphism of exact sequences 

$$
\begin{array}{ccccccc}
1\longrightarrow & \T(S)&\longrightarrow &\mme (S) &\longrightarrow
&\spe (2g,\Z)&\longrightarrow 1\\
& & & & & & \\
&\downarrow \,\approx& & \downarrow \, \approx & & \downarrow \, \approx & \\
& & & & & & \\
1\longrightarrow &\Inn(\T(S)))&\longrightarrow &\Aut(\T(S))&\longrightarrow
&\Out(\T(S))&\longrightarrow 1,
\end{array}$$
where the first exact sequence is an obvious version of the exact sequence (\ref{equation:basic})
and the second exact sequence is the usual one for a centerless group, with $\Inn$ denoting the group of inner automorphisms induced by conjugation.
\end{theorem}

Note that the conclusion of Theorem 1 is false when $\genus(S)=2$, as Mess
\cite{Me} has shown that in this case $\T(S)$ is a countably generated free
group. Theorem \ref{theorem:autos} was inspired by the theorem of 
Ivanov \cite{i0}, \cite{i-aut} (see also \cite{mc-aut}) to the effect 
that $\Out(\mme(S))=1$ for $\genus(S)\geq 3$.

\bigskip
\noindent
{\bf The commensurator. }Recall that the {\em (abstract) commensurator
group} $\Comm(\Gamma)$ of a group $\Gamma$ is defined to be the set of
equivalence classes of isomorphisms $\phi:H\rightarrow N$ between finite
index subgroups $H,N$ of $\Gamma$, where $\phi_1:H_1\rightarrow N_1$ 
is equivalent to $\phi_2:H_2\rightarrow N_2$ if $\phi_1=\phi_2$ on
$H_1\cap H_2$. The composition of homomorphisms induces a natural multiplication
on $\Comm(\Gamma)$, which turns $\Comm(\Gamma)$ into a group.  

The group $\Comm(\Gamma)$ 
is in general much larger than $\Aut(\Gamma)$; for example
$\Aut(\Z^n)=\GL(n,\Z)$ whereas $\Comm(\Z^n)=\GL(n,\Q)$.  
The group $\Comm(\Gamma)$ was computed
for mapping class groups of surfaces by Ivanov \cite{Iv2}.  
For the Torelli group we have the following:

\begin{theorem}[Commensurator Theorem]
\label{theorem:comm}
Let $S$ be a closed surface of $\genus(S)\geq 5$.  
Then the natural injection 
$$\mme (S)\to \Comm(\T_S)$$
is an isomorphism.
\end{theorem}

\bigskip
\noindent
{\bf The non-arithmeticity. }Our proof of the non-arithmeticity of
$\T_S$ is inspired by the proof of the non-arithmeticity of $\Mod (S)$
given in \cite{Iv2}. Theorem \ref{theorem:comm} implies that $\T
_S$ is normal in its abstract commensurator. Using this fact, it is easy
to see that $\T_S$ cannot be arithmetic. As an illustrating example,
notice that the arithmetic group $\GL(n,\Z)$ is not normal in its
abstract commensurator $\Comm(\GL(n,\Z))=\GL(n,\Q)$.

\bigskip
\noindent
{\bf The idea of proof of Theorems \ref{theorem:autos} and \ref{theorem:comm}. }The main step in the proof of Theorems \ref{theorem:autos} and \ref{theorem:comm} is to prove 
that the given isomorphism $\Phi: \Gamma_1 \to \Gamma_2$ between two subgroups of finite index in $\T_S$ induces an automorphism of the 
Torelli Geometry $\TG(S)$; Theorem \ref{automorphisms of geometry} can then be
applied to produce a mapping class inducing $\Phi$.  
First of all, we note that the isotopy
class of a curve is uniquely determined by (any non-zero power of) the Dehn twist about that
curve. If one has a purely algebraic characterization of (the powers of) Dehn twists 
in the Torelli group, this can be used to produce an automorphism of $\TG(S)$.

We will denote the Dehn twist about a curve $\gamma$ by $T_\gamma$, and define
the Dehn twist about a bounding pair $\{a,b\}$ as
$T_{ab^{-1}}=T_a T_{b^{-1}}$. Dehn twist about a separating curve or
bounding pair will be called a {\em simple twist}. For a group $\Gamma$
and element $f\in \Gamma$, let $C_\Gamma(f)$ denote the subgroup of
elements $g\in \Gamma$ commuting with $f$, and let $Z(\Gamma)$ denote
the center of $\Gamma$. The following proposition gives us an algebraic 
characterization of (powers of) simple twists.

\begin{proposition}[Characterizing simple twists in \boldmath$\T(S)$]
\label{proposition:characterizing:twists}
Let $S$ be a closed surface with $\genus(S)\geq 3$, 
and let $f\in \T(S)$ be nontrivial. Then $f$ is a power of
a simple twist if and only if both of the following hold:
\begin{enumerate}
\item $Z(C_{\T(S)}(f))=\Z$, and
\item The maximum of ranks of abelian subgroups of $\T(S)$ that contain
$f$ is $2g-3$.
\end{enumerate}
\end{proposition}

\protect From Proposition \ref{proposition:characterizing:twists} 
one can easily deduce that any isomorphism $\Phi:\Gamma_1\to \Gamma_2$ between two subgroups $\Gamma_1$, $\Gamma_2$ of finite index in $\T _S$
takes powers of simple twists to powers of simple twists. The main tool in proving Proposition
\ref{proposition:characterizing:twists} is the Thurston normal form
theory for mapping classes extended to abelian subgroups of $\Mod(S)$, for example, in \cite{Iv1}.

Next, we need further to distinguish 
purely algebraically between the two types of simple twists. 

\begin{proposition}[Characterizing bounding pairs]
\label{proposition:characterizing:boundingpairs}
Let $g\geq 4$, and let $R$ be the set of powers of 
simple twists in $\T(S)$.  Then the following are equivalent.
\begin{enumerate}
\item $f\in R$ is a power of a Dehn twist about a bounding pair.
\item There exist $g,h\in R$, distinct from $f$ and from each other, so
that the group generated by $f,g,h$ is isomorphic to $\Z^2$.
\end{enumerate}
\end{proposition}

Finally, we need to give a purely algebraic characterization of marked triangles in $\TG (S)$.

\begin{proposition}[Characterizing marked triangles]
\label{proposition:characterizing:markedtriangles}
Let $g\geq 4$, and let $R$ be the set of powers of 
simple twists in $\T(S)$.  Then the following are equivalent.
\begin{enumerate}
\item $f,g,h\in R$ are powers of  Dehn twists about bounding pairs $\beta_1$, $\beta_2$,
$\beta_3$ such that $\{\beta_1,\beta_2,\beta_3\}$ is a marked triangle.
\item The groups generated by elements $f,g,h$ and by any two of them are all isomorphic to $\Z^2$.
\end{enumerate}
\end{proposition}

The proofs of the last two Propositions are similar, and are based on an analysis of all possible
configurations of three separating circles and bounding pairs on a surface.

\noindent
Benson Farb:\\
Dept. of Mathematics, University of Chicago\\
5734 University Ave.\\
Chicago, Il 60637\\
E-mail: farb@math.uchicago.edu

\medskip
\noindent
Nikolai V. Ivanov:\\
Michigan State University\\
Dept. of Math, Wells Hall\\
East Lansing, MI 48824-1027\\
E-mail: ivanov@math.msu.edu
\end{document}